\documentclass[11pt,a4paper]{article}
\usepackage[T1]{fontenc}
\usepackage[margin=28mm]{geometry}
\usepackage{enumitem}
\usepackage{amsmath,amssymb,amsthm,booktabs,array,needspace}
\usepackage[colorlinks=true,linkcolor=blue,citecolor=blue,urlcolor=blue]{hyperref}
\hypersetup{pdftitle={Extremal spectral results of outerplanar graphs without starlike trees},pdfauthor={Fulong Ye, Bofeng Huo, Yuxiang Liu and Ligong Wang},pdfsubject={Spectral extremal graph theory}}
\IfFileExists{glyphtounicode.tex}{\input{glyphtounicode}\pdfgentounicode=1}{}
\newtheorem{theorem}{Theorem}[section]
\newtheorem{lemma}[theorem]{Lemma}
\newtheorem{corollary}[theorem]{Corollary}

\theoremstyle{definition}
\newtheorem{definition}[theorem]{Definition}
\newtheorem{problem}{Problem}
\theoremstyle{remark}

\newcommand{\spex}{\operatorname{spex}_{\mathcal{OP}}}
\newcommand{\SPEX}{\operatorname{SPEX}_{\mathcal{OP}}}
\newcommand{\spexP}{\operatorname{spex}_{\mathcal P}}
\newcommand{\SPEXP}{\operatorname{SPEX}_{\mathcal P}}
\newcommand{\HH}{H_{\mathcal{OP}}}
\newcommand{\tH}{\widetilde H_{\mathcal{OP}}}
\newcommand{\floor}[1]{\left\lfloor #1\right\rfloor}

\newcommand{\claim}[1]{\par\Needspace{3\baselineskip}\medskip\noindent\textbf{Claim #1.}\ }
\newcommand{\case}[1]{\par\Needspace{3\baselineskip}\medskip\noindent\textbf{Case #1.}\ }

\title{Extremal spectral results of outerplanar graphs without starlike trees\thanks{Supported by the National Natural Science Foundation of China (Nos. 12271439 and 11961055).}}
\author{Fulong Ye$^{a,b,c}$, Bofeng Huo$^{d}$,
Yuxiang Liu$^{a,c}$, Ligong Wang$^{a,c}$%
\thanks{Corresponding author.}\\[1.2ex]
\small $^{a}$School of Mathematics and Statistics,
Northwestern Polytechnical University,\\
\small Xi'an, Shaanxi 710129, P.R. China.\\[0.5ex]
\small $^{b}$School of Mathematics and Statistics,
Qinghai Minzu University,\\
\small Xining, Qinghai 810007, P.R. China.\\[0.5ex]
\small $^{c}$Xi'an--Budapest Joint Research Center
for Combinatorics,\\
\small Northwestern Polytechnical University,
Xi'an, Shaanxi 710129, P.R. China.\\[0.5ex]
\small $^{d}$College of Mathematics and Statistics,
Qinghai Normal University,\\
\small Xining, Qinghai 810008, P.R. China.\\[0.5ex]
\small E-mail: \texttt{yefl\_qhmu@163.com},
\texttt{huobofeng@qhnu.edu.cn},\\
\small \texttt{yxliumath@163.com},
\texttt{lgwangmath@163.com}}
\date{}
\begin{document}
\maketitle
\vspace{-8mm}
\begin{abstract}
A graph $G$ is $F$-free if it does not contain $F$ as a subgraph.
Let $\operatorname{SPEX}_{\mathcal{OP}}(n,F)$ denote the set of
$n$-vertex $F$-free outerplanar graphs with maximum spectral radius.
Define $S(l_1,l_2,\ldots,l_t)$ as the starlike tree obtained from
$t$ vertex-disjoint paths of lengths $l_1,l_2,\ldots,l_t$ by
identifying one endpoint of each path, and $C_{l_1,l_2,\ldots,l_t}$
as the graph obtained from $t$ vertex-disjoint cycles of lengths
$l_1,l_2,\ldots,l_t$ by identifying one vertex from each cycle.
When $l_1=l_2=\cdots=l_t=l$, we abbreviate these graphs as
$S(t\cdot l)$ and $C_{t\cdot l}$, respectively.
Previous results determined
$\operatorname{SPEX}_{\mathcal{OP}}(n,S(t\cdot l))$
for $t\ge4$ [Graphs Combin. 42 (2026) 80]
and for $t=3$, in each case for sufficiently large $n$.
For $C_{l_1,l_2,\ldots,l_t}$, the corresponding extremal graphs
are known for $t=2$ and sufficiently large $n$.
The sets $\operatorname{SPEX}_{\mathcal{OP}}(n,C_{t\cdot l})$
were also determined for sufficiently large $n$
[Discrete Math. 349 (2026) 114673].
In this paper, we determine
$\operatorname{SPEX}_{\mathcal{OP}}(n,S(l_1,l_2,\ldots,l_t))$
for sufficiently large $n$ when the branch lengths are not all equal.
Our result for $t\ge4$ also includes equal branch lengths.
As applications, we determine
$\operatorname{SPEX}_{\mathcal{OP}}(n,C_{l_1,l_2,\ldots,l_t})$
for arbitrary cycle lengths and sufficiently large $n$.
We also determine
$\operatorname{SPEX}_{\mathcal{OP}}(n,F)$ for certain linear
forests $F=\bigcup_{i=1}^tP_{l_i}$ and sufficiently large $n$.
Finally, we apply our main theorems to obtain the known
results on outerplanar graphs forbidding matchings or
theta graphs for sufficiently large $n$.
\end{abstract}

\section{Introduction}
All graphs considered here are finite and undirected without
loops or multiple edges.
For a graph $G=(V(G),E(G))$, let $V(G)$ and $E(G)$ be its
\emph{vertex set} and \emph{edge set}, respectively.
The \emph{order} of $G$ is $n=|V(G)|$, and its \emph{size}
is $e(G)=|E(G)|$.
The \emph{adjacency matrix} of $G$ is the $n\times n$ matrix
$A(G)=(a_{uv})_{u,v\in V(G)}$, whose rows and columns are
indexed by the vertices of $G$.
Its entries satisfy $a_{uv}=1$ if $uv\in E(G)$ and
$a_{uv}=0$ if $uv\notin E(G)$.
Since $A(G)$ is real and symmetric, all its eigenvalues are real.
The largest eigenvalue of $A(G)$, denoted by $\rho(G)$, is called the
\emph{spectral radius} of $G$.

For $v\in V(G)$, its \emph{neighborhood} $N_G(v)$ is the set
of vertices adjacent to $v$, and its \emph{degree} is
$d_G(v)=|N_G(v)|$.
For $U\subseteq V(G)$, let $G[U]$ denote the
\emph{subgraph induced by $U$}, consisting of
the vertices in $U$ and all edges of $G$ with
both endpoints in $U$.
Let $G-U=G[V(G)\setminus U]$ denote the graph
obtained from $G$ by deleting the vertices in $U$
and all edges incident with them.
For $U=\{v\}$, we use $G-v$.
A graph $H$ is a \emph{subgraph} of $G$,
denoted by $H\subseteq G$, if
$V(H)\subseteq V(G)$ and $E(H)\subseteq E(G)$.
Let $P_s$, $C_s$ and $K_s$ denote the \emph{path},
\emph{cycle} and \emph{complete graph} on $s$ vertices,
respectively.
Let $K_{p,q}$ denote the \emph{complete bipartite graph}
with two parts of sizes $p$ and $q$.
For two vertex-disjoint graphs $G$ and $H$, their \emph{union}
is denoted by $G\cup H$.
Their \emph{join}, denoted by $G\vee H$, is obtained by
adding all edges between $V(G)$ and $V(H)$ to $G\cup H$.
For a positive integer $q$, let $qG$ denote the disjoint union
of $q$ copies of $G$.

For later use, we introduce the graphs considered in this paper.
A \emph{linear forest} is a graph whose components
are paths.
For integers $t\ge3$ and $l_1\ge l_2\ge \cdots\ge l_t\ge1$,
let $S(l_1,l_2,\ldots,l_t)$ be the starlike tree obtained
by identifying one endpoint of each of $t$ otherwise disjoint paths of lengths $l_1,l_2,\ldots,l_t$.
When $l_1=l_2=\cdots=l_t=l$, we denote this tree by $S(t\cdot l)$.
This tree has $1+\sum_{i=1}^t l_i$ vertices. Its unique vertex of degree $t$ is called the center, say $w$, and
$
S(l_1,l_2,\ldots,l_t)-w\cong\bigcup_{i=1}^tP_{l_i}.
$
For integers $t\ge2$ and
$l_1\ge l_2\ge\cdots\ge l_t\ge1$, let
$F(l_1,l_2,\ldots,l_t)=\bigcup_{i=1}^tP_{l_i}$.
For integers $t\ge2$ and $l_1\ge l_2\ge\cdots\ge l_t\ge3$,
let $C_{l_1,l_2,\ldots,l_t}$ be the graph obtained by identifying one vertex from each of $t$ vertex-disjoint cycles of lengths $l_1,l_2,\ldots,l_t$.
When $l_1=l_2=\cdots=l_t=l$, we denote this graph by $C_{t\cdot l}$.
For integers $r,p,q$ with $1\le r\le p\le q$ and $p\ge2$,
let $\theta_{r,p,q}$ be the graph formed by three internally
vertex-disjoint paths of lengths $r,p,q$ with the same endpoints.
It has $r+p+q-1$ vertices.
In particular, $\theta_{1,p,q}$ is a cycle of length $p+q$
with an added chord.
For an integer $k\ge4$, let
$\Theta_k=\{\theta_{1,p,q}:2\le p\le q,\ p+q=k\}$.

For a family $\mathcal F$ of graphs, a graph $G$ is called
$\mathcal F$-free if it does not contain any member of
$\mathcal F$ as a subgraph.
When $\mathcal F=\{F\}$, we simply say that
$G$ is $F$-free.
The well-known Brualdi--Solheid problem~\cite{BS} seeks to
determine the maximum spectral radius within a prescribed
graph class and characterize the corresponding extremal graphs.
In 2010, Nikiforov~\cite{Npaths} proposed a forbidden-subgraph version of this problem,
asking for the maximum spectral radius of an $F$-free graph of order $n$ for a given graph $F$.
During the past decades, this problem has been extensively studied for various forbidden subgraphs,
including complete graphs~\cite{BN},  paths and  cycles~\cite{CDTeven,Npaths},
star forests~\cite{CLZstar}, linear forests~\cite{CLZforest},
and trees~\cite{CDTtrees,FLSZ,LBWtrees}.
Further results concern forbidden spanning linear forests~\cite{ZLPW} and disjoint unions of three-vertex paths~\cite{ZPW}.
Spectral extremal problems with a fixed number of edges have also been studied in~\cite{LZZsize,LLLYtheta,LWsize,LLZmantel,ZYWstar}.
For more results, we refer readers to the surveys~\cite{CZsurvey,LLFsurvey,Nsurvey}.

A graph is planar if it can be drawn in the plane without edge crossings.
A planar graph is outerplanar
if it has a plane drawing in which all vertices lie on the boundary of the outer face.
Let $\spexP(n,F)$ and $\spex(n,F)$
denote the maximum spectral radius over all $n$-vertex $F$-free planar and outerplanar graphs,
respectively,
and let $\SPEXP(n,F)$ and $\SPEX(n,F)$ denote the corresponding sets of extremal graphs.

Spectral extremal problems for planar graphs date back to Schwenk and Wilson~\cite{SW}.
Cvetkovi\'c and Rowlinson~\cite{CR} conjectured that $K_1\vee P_{n-1}$ maximizes the spectral radius of outerplanar graphs,
while Boots and Royle~\cite{BR},
and independently Cao and Vince~\cite{CV},
conjectured that $K_2\vee P_{n-2}$ is the unique
$n$-vertex planar graph with maximum spectral
radius for $n\ge9$.
Tait and Tobin~\cite{TT} proved both conjectures for sufficiently large $n$.
Subsequently, Lin and Ning~\cite{LN} completely solved the outerplanar problem,
and Liu, Ning and Wang~\cite{LNW} completely solved the planar problem.

These results motivated spectral extremal problems for planar
and outerplanar graphs with forbidden subgraphs.
For planar graphs, results on $\SPEXP(n,F)$ have been obtained
for $F\in\{tC_l,S(t\cdot l),C_{t\cdot l},C_{l_1,l_2}\}$
and the families $t\mathcal C$ and $\Theta_k$;
see~\cite{FLS,LLYC,YL,YLM,ZW},
where $t\mathcal C$ is the family of disjoint unions of $t$ cycles.
Zhai and Liu~\cite{ZLedgecycles} also determined the unique
spectral extremal planar graph without $k$ edge-disjoint cycles.
For outerplanar graphs, $\SPEX(n,F)$ has been determined for
$F\in\{S(t\cdot l),C_{t\cdot l},C_{l_1,l_2},
(k+1)K_2,\theta_{1,p,q}\}$;
see~\cite{LLYC,WSS,YL,YLM}.
These results hold for sufficiently large $n$ under the
respective parameter conditions.
We introduce some notation and state the relevant problems
and known results below.
\begin{definition}[\cite{YL}]\label{def:H}
For positive integers $n,n_1,n_2$ with $n-1\ge n_1\ge n_2\ge1$, define
\[
\HH(n_1,n_2)=
\begin{cases}
P_{n_1}\cup\dfrac{n-1-n_1}{n_2}P_{n_2},
&\text{if }n_2\mid(n-1-n_1),\\[7pt]
P_{n_1}\cup\floor{\dfrac{n-1-n_1}{n_2}}P_{n_2}
\cup P_{\,n-1-n_1-\floor{\frac{n-1-n_1}{n_2}}n_2},
&\text{otherwise}.
\end{cases}
\]
\end{definition}

\begin{definition}[\cite{YLW}]\label{def:Htilde}
For positive integers $n,n_1,n_2$ with $n-1\ge2n_1$ and $n_1\ge n_2\ge1$, define
\[
\tH(n_1,n_2)=
\begin{cases}
2P_{n_1}\cup\dfrac{n-1-2n_1}{n_2}P_{n_2},
&\text{if }n_2\mid(n-1-2n_1),\\[7pt]
2P_{n_1}\cup\floor{\dfrac{n-1-2n_1}{n_2}}P_{n_2}
\cup P_{\,n-1-2n_1-\floor{\frac{n-1-2n_1}{n_2}}n_2},
&\text{otherwise}.
\end{cases}
\]
\end{definition}
Yin, Li and Meng~\cite{YLM} posed the following spectral extremal problem for $S(l_1,l_2,\ldots,l_t)$.
\begin{problem}[\cite{YLM}]\label{prob:starlike}
Characterize the planar and outerplanar spectral extremal graphs of order $n$ without a starlike tree $S(l_1,l_2,\ldots,l_t)$ for sufficiently large $n$.
\end{problem}
For planar and outerplanar graphs forbidding
$C_{l_1,l_2,\ldots,l_t}$, we pose the following
analogous spectral extremal problem.
\begin{problem}\label{prob:cycles}
Characterize the planar and outerplanar spectral extremal graphs of order $n$ without $C_{l_1,l_2,\ldots,l_t}$ for sufficiently large $n$.
\end{problem}

We first state the known results relevant to
Problem~\ref{prob:starlike}.
For starlike trees whose branches all have the same
length, Yin, Li and Meng~\cite{YLM} solved the case $t\ge4$,
and Ye, Liu and Wang~\cite{YLW} solved the case $t=3$.
\begin{theorem}[\cite{YLM}]\label{thm:knownstarlike}
For integers $n,t,l$ with $t\ge4$, $l\ge2$ and
$n\ge\max\{1.27\times10^7,6.5025\times2^{tl+1}\}$,
let $G$ be an $n$-vertex $S(t\cdot l)$-free outerplanar graph.
Then $\rho(G)\le \rho\bigl(K_1\vee\HH(tl-1,l-1)\bigr),$ with equality if and only if $G\cong K_1\vee\HH(tl-1,l-1)$.
\end{theorem}
\begin{theorem}[\cite{YLW}]\label{thm:equal}
For integers $n,l$ with $l\ge2$ and
$n\ge\max\{1.27\times10^7,6.5025\times2^{3l+1},(\sqrt{65+3l}+8)^2\}$,
let $G$ be an $n$-vertex $S(3\cdot l)$-free outerplanar graph.
Then $\rho(G)\le \rho\bigl(K_1\vee\tH(2l-1,l-1)\bigr),$
with equality if and only if $G\cong K_1\vee\tH(2l-1,l-1)$.
\end{theorem}

For Problem~\ref{prob:cycles}, Yin and Li~\cite{YL} determined
the extremal graphs when $l_1=l_2=\cdots=l_t$,
and Li et al.~\cite{LLYC} determined them when
$t=2$ and $l_1>l_2$.
\begin{theorem}[\cite{YL}]\label{thm:knowncycles}
For integers $n,t,l$ with $t\ge2$, $l\ge3$ and
$n\ge\max\{1.27\times10^7,6.5025\times2^{t(l-1)+1}\}$,
let $G$ be an $n$-vertex $C_{t\cdot l}$-free outerplanar graph.
Then $\rho(G)\le \rho\bigl(K_1\vee\HH(tl-t-1,l-2)\bigr),$
with equality if and only if $G\cong K_1\vee\HH(tl-t-1,l-2)$.
\end{theorem}
\begin{theorem}[\cite{LLYC}]\label{thm:twocycles}
For integers $n,l_1,l_2$ with $l_1>l_2\ge3$ and
$n>\max\{400,\allowbreak\frac54(l_1+l_2-1),\allowbreak
5\times2^{l_1}+1,\allowbreak
\frac{(l_2-2)(2l_1^2-14l_1+3l_1l_2-6l_2+19)}{l_1-l_2}\}$,
let $G$ be an $n$-vertex $C_{l_1,l_2}$-free outerplanar graph.
Then $\rho(G)\le
\rho\bigl(K_1\vee\HH(l_1-2,l_1-2)\bigr),$
with equality if and only if
$G\cong K_1\vee\HH(l_1-2,l_1-2).$
\end{theorem}

Motivated by Problems~\ref{prob:starlike}
and~\ref{prob:cycles} and the preceding results,
we characterize the $n$-vertex  $S(l_1,l_2,\ldots,l_t)$-free outerplanar
graphs with maximum spectral radius for every
remaining nonstar starlike tree  $S(l_1,l_2,\ldots,l_t)$.
As an application, we characterize the $n$-vertex
$C_{l_1,l_2,\ldots,l_t}$-free outerplanar graphs with maximum
spectral radius for $t\ge2$ and $l_1\ge l_2\ge\cdots\ge l_t\ge3$.
For an integer $L\ge4$, let
\begin{equation}\label{eq:n0}
n_0=\max\{1.27\times10^7,\ 6.5025\times2^{L+1},\ (8L^2+8)^2\}.
\end{equation}
This common sufficient bound $n_0$ simplifies the proofs.
It is not claimed to be optimal.

Our first theorem determines
$\SPEX(n,S(l_1,l_2,l_3))$
when the three branch lengths are not all equal.
\begin{theorem}\label{thm:three}
For integers $l_1,l_2,l_3$ with
$l_1\ge l_2\ge l_3\ge1$ and $l_1>l_3$, let
$L=l_1+l_2+l_3$.
For an integer $n\ge n_0$ \textup{(see~\eqref{eq:n0})}, let $G$ be an $n$-vertex
$S(l_1,l_2,l_3)$-free outerplanar graph. Then $\rho(G)\le \rho\bigl(K_1\vee H^*\bigr),$
where
\[
H^*=
\begin{cases}
\HH(l_1-1,l_1-1),
& \text{if } l_1>l_2 \text{ and } l_1\le l_2+l_3+1,\\
\HH(l_1-1,l_1-2),
& \text{if } l_1=l_2+l_3+2,\\
\HH(l_1-2,l_1-2),
& \text{if } l_1\ge l_2+l_3+3,\\
\HH(l_1+l_3,l_1-1),
& \text{if } l_1=l_2>l_3.
\end{cases}
\]
Moreover, equality holds if and only if $G\cong K_1\vee H^*$.
\end{theorem}

For a nonstar starlike tree $S(l_1,l_2,\ldots,l_t)$ with $t\ge4$,
the unique $n$-vertex $S(l_1,l_2,\ldots,l_t)$-free outerplanar
graph with maximum spectral radius depends only on $l_1$
and the number of branches of length $l_1$,
provided that $n$ is sufficiently large.
\begin{theorem}\label{thm:many}
For integers $t,l_1,l_2,\ldots,l_t$ with $t\ge4$,
$l_1\ge l_2\ge\cdots\ge l_t\ge1$ and $l_1\ge2$, let
$L=\sum_{i=1}^t l_i$ and $h=|\{i:l_i=l_1,\ i\in\{1,2,\ldots,t\}\}|$.
For an integer $n\ge n_0$ \textup{(see~\eqref{eq:n0})}, let $G$ be an $n$-vertex
$S(l_1,l_2,\ldots,l_t)$-free outerplanar graph. Then
$\rho(G)\le
\rho\bigl(K_1\vee\HH(hl_1-1,l_1-1)\bigr),$
with equality if and only if
$G\cong K_1\vee\HH(hl_1-1,l_1-1).$
\end{theorem}

When $l_1=l_2=\cdots=l_t$, we have $h=t$.
By Theorem~\ref{thm:many},
the unique extremal graph is $K_1\vee\HH(tl_1-1,l_1-1)$,
which agrees with Theorem~\ref{thm:knownstarlike}
for starlike trees whose branches all have the
same length.
Hence Theorem~\ref{thm:many} applies whether or not
the branch lengths are all equal when $t\ge4$.
 Combining Theorems~\ref{thm:equal},
\ref{thm:three}, and~\ref{thm:many}, we summarize all cases
in the following table, where the unique extremal graph
is $K_1\vee H^*$ for sufficiently large $n$.

\begin{table}[ht]
\centering
\caption{Extremal forests for nonstar starlike trees.}\label{tab:allstarlike}\smallskip
\begin{tabular}{@{}lll@{}}
\toprule
Number of branches & Branch lengths & $H^*$\\
\midrule
$t=3$ & $l_1=l_2=l_3$ & $\tH(2l_1-1,l_1-1)$\\
$t=3$ & $l_1>l_2$, $l_1\le l_2+l_3+1$ & $\HH(l_1-1,l_1-1)$\\
$t=3$ & $l_1=l_2+l_3+2$ & $\HH(l_1-1,l_1-2)$\\
$t=3$ & $l_1\ge l_2+l_3+3$ & $\HH(l_1-2,l_1-2)$\\
$t=3$ & $l_1=l_2>l_3$ & $\HH(l_1+l_3,l_1-1)$\\ $t\ge4$ & $l_1\ge l_2\ge\cdots\ge l_t\ge1$, $l_1\ge2$,
  $h=|\{i:l_i=l_1\}|$
& $\HH(hl_1-1,l_1-1)$\\
\bottomrule
\end{tabular}
\end{table}

Let $H$ be a linear forest and let $u$ be the vertex
corresponding to $K_1$ in $K_1\vee H$.
The case $t=3$ requires several subcases because the center
of $S(l_1,l_2,l_3)$ may be either $u$ or a vertex of $H$.
When $t\ge4$, the center must be $u$, since every vertex
of $H$ has degree at most $3$ in $K_1\vee H$.

As an application of Theorem~\ref{thm:many},
we also solve the outerplanar part of Problem~\ref{prob:cycles}.
The result includes Theorems~\ref{thm:knowncycles} and~\ref{thm:twocycles} as special cases.
\begin{theorem}\label{thm:bouquet}
For integers $t,l_1,l_2,\ldots,l_t$ with $t\ge2$ and
$l_1\ge l_2\ge\cdots\ge l_t\ge3$, let
$L=\sum_{i=1}^t(l_i-1)+2$ and $h=|\{i:l_i=l_1,\ i\in\{1,2,\ldots,t\}\}|$.
For an integer $n\ge n_0$ \textup{(see~\eqref{eq:n0})}, let $G$ be an $n$-vertex
$C_{l_1,l_2,\ldots,l_t}$-free outerplanar graph. Then
$
\rho(G)\le
\rho\bigl(
K_1\vee\HH(h(l_1-1)-1,l_1-2)
\bigr),$
with equality if and only if
$G\cong
K_1\vee\HH(h(l_1-1)-1,l_1-2).
$
\end{theorem}

Yin, Li and Meng~\cite{YLM} also posed the following problem for linear forests as Problem~6 in their paper.
\begin{problem}[\cite{YLM}]\label{prob:forest}
Characterize the planar and outerplanar spectral extremal graphs of order $n$
without a  linear forest $F$ for sufficiently large $n$.
\end{problem}

Yin, Li and Meng~\cite{YLM} determined
$\SPEX(n,tP_l)$ for sufficiently large $n$,
where $tP_l$ is the disjoint union of $t\ge2$ copies of $P_l$
with $l\ge3$.
By constructing starlike trees containing $F(l_1,l_2,\ldots,l_t)$
and applying Theorems~\ref{thm:three} and~\ref{thm:many},
we obtain the following partial answer to Problem~\ref{prob:forest} for outerplanar graphs.

\begin{theorem}\label{thm:forest}
For integers $t,l_1,l_2,\ldots,l_t$ with $t\ge2$,
$l_1\ge l_2\ge\cdots\ge l_t\ge2$ and either
$l_1\le2l_2$ or $l_1$ is even, let
$L=\sum_{i=1}^t l_i$ and $h=|\{i:l_i=l_2,\ i\in\{1,2,\ldots,t\}\}|$.
For an integer $n\ge\max\{1.27\times10^7,
6.5025\times2^{L+2},(8(L+1)^2+8)^2\}$,
let $G$ be an $n$-vertex
$F(l_1,l_2,\ldots,l_t)$-free outerplanar graph. Then $\rho(G)\le\rho(K_1\vee H^*)$, where
\[
H^*=
\begin{cases}
H_{\mathcal{OP}}((h-1)l_1-1,l_1-1),
& \text{if } l_1=l_2,\\[3pt]
H_{\mathcal{OP}}(hl_2-1,l_2-1),
& \text{if } l_2<l_1\le2l_2-1,\\[3pt]
H_{\mathcal{OP}}((h+1)l_2-1,l_2-1),
& \text{if } l_1=2l_2,\\[3pt]
H_{\mathcal{OP}}\left(
\frac{l_1}{2}-1,\frac{l_1}{2}-1\right),
& \text{if } l_1>2l_2 \text{ and } l_1 \text{ is even},
\end{cases}
\]
with equality  if and only if $G\cong K_1\vee H^*$.
\end{theorem}

An outline of the paper is as follows.
Section~2 presents some preliminary results.
Sections~3 and~4 prove Theorems~\ref{thm:three} and~\ref{thm:many}, respectively.
Section~5 proves Theorems~\ref{thm:bouquet} and~\ref{thm:forest}.
\section{Preliminaries}
One of the fundamental properties of outerplanar graphs is that
they contain no $K_{2,3}$ minor. Consequently, any two distinct
vertices of an outerplanar graph have at most two common neighbors,
since three common neighbors together with these two vertices would
contain  $K_{2,3}$ as a subgraph. Moreover, every outerplanar graph of order
$n\ge2$ has at most $2n-3$ edges. We refer the reader to~\cite{BM}
for more results on outerplanar graphs. Finally, if $H$ is a linear
forest, then $K_1\vee H$ is outerplanar.

For a graph $G$, let $\boldsymbol{x}=(x_1,x_2,\ldots,x_n)^{\mathrm T}$ be an eigenvector of
$A(G)$ corresponding to $\rho(G)$. For each vertex $v\in V(G)$, the
eigenvalue equation is
\[
\rho(G)x_v=\sum_{u\in N_G(v)}x_u.
\]
By the Rayleigh principle, we have
\[
\rho(G)
=\max_{\boldsymbol{x}\in\mathbb{R}^n\setminus\{\boldsymbol{0}\}}
\frac{\boldsymbol{x}^{\mathrm T}A(G)\boldsymbol{x}}
     {\boldsymbol{x}^{\mathrm T}\boldsymbol{x}}
=\max_{\boldsymbol{x}\in\mathbb{R}^n\setminus\{\boldsymbol{0}\}}
\frac{2\displaystyle\sum_{uv\in E(G)}x_ux_v}
     {\boldsymbol{x}^{\mathrm T}\boldsymbol{x}}.
\]
If $G$ is connected, the Perron--Frobenius theorem guarantees a
positive eigenvector $\boldsymbol{x}=(x_1,x_2,\ldots,x_n)^{\mathrm T}$ corresponding to $\rho(G)$.
If $G'$ has the same vertex set as $G$, then
\[
 \rho(G')-\rho(G)\ge
 \frac{2}{\boldsymbol{x}^{\mathrm T}\boldsymbol{x}}
 \left(
 \sum_{vw\in E(G')\setminus E(G)}x_vx_w-
 \sum_{vw\in E(G)\setminus E(G')}x_vx_w
 \right).
\]

\begin{lemma}[{\cite{LiFeng}}]\label{lem:proper-spanning}
Let $G$ be a connected graph and let $H$ be a proper spanning
subgraph of $G$. Then $\rho(H)<\rho(G)$.
\end{lemma}

\begin{lemma}[\cite{YL}]\label{lem:YLstructure}
Let $F$ be an outerplanar graph contained in
$K_1\vee P_{n-1}$ but not contained in $K_1\vee\bigl((t-1)K_2\cup(n-2t+1)K_1\bigr),$
where $1\le t\le \frac{n-1}{2}$ and
$n\ge
\max\left\{
1.27\times10^7,\,
\left(\sqrt{64+|V(F)|}+8\right)^2
\right\}.$
Let $G\in\SPEX(n,F)$ be connected, and let
$\boldsymbol{x}$ be a positive eigenvector
corresponding to $\rho(G)$ with
$
\max_{v\in V(G)}x_v=1.
$
Then there exists a vertex $u\in V(G)$
such that $N_G(u)=V(G)\setminus\{u\}$ and $x_u=1.$
In particular,
\begin{enumerate}
[
  label=(\roman*),
  leftmargin=*,
  topsep=1pt,
  itemsep=3pt,
  parsep=0pt,
  partopsep=0pt
]
\item $G$ contains a copy of $K_{1,n-1}$.
\item The subgraph $G[N_G(u)]$ is a linear forest.
\end{enumerate}
\end{lemma}

\begin{definition}[\cite{YL}]\label{def:transformation}
For integers $s_1,s_2$ with $s_1\ge s_2\ge1$, let
$H=P_{s_1}\cup P_{s_2}\cup H_0$, where $H_0$ is a linear forest.
We say that $H'$ is an $(s_1,s_2)$-transformation of $H$ if
\[
H'=
\begin{cases}
P_{s_1+1}\cup P_{s_2-1}\cup H_0,
& \text{if }s_2\ge2,\\
P_{s_1+1}\cup H_0,
& \text{if }s_2=1.
\end{cases}
\]
\end{definition}
Clearly, $H'$ is a linear forest. Hence $K_1\vee H'$ is
outerplanar.

\begin{lemma}[\cite{YL}]\label{lem:transfer}
For integers $n,s_1,s_2$ with $s_1\ge s_2\ge1$ and
$n\ge\max\{1.27\times10^7,6.5025\times2^{s_2+2}\}$,
let $H$ be a linear forest of order $n-1$ and let $H'$
be an $(s_1,s_2)$-transformation of $H$. Then
$\rho(K_1\vee H')>\rho(K_1\vee H).$
\end{lemma}

\begin{lemma}\label{lem:reduction}
For integers $t,l_1,l_2,\ldots,l_t$ with $t\ge3$,
$l_1\ge l_2\ge\cdots\ge l_t\ge1$ and $l_1\ge2$, let
$F=S(l_1,l_2,\ldots,l_t)$ and $L=\sum_{i=1}^t l_i$.
For an integer $n\ge n_0$ \textup{(see~\eqref{eq:n0})}, let $G\in\SPEX(n,F)$.
Then $G$ is connected
and has the form $G=K_1\vee H$, where $H$ is a linear forest.
\end{lemma}
\begin{proof}
We first show that $G$ is connected.
Since $F$ is not a star, $K_{1,n-1}$ is an $F$-free
outerplanar graph. By the extremality of $G$, we obtain
$\rho:=\rho(G)\ge\rho(K_{1,n-1})=\sqrt{n-1}$.
Suppose that $G$ is disconnected.
Let $G_1$ be a connected component with $\rho(G_1)=\rho$.
Choose a positive Perron vector on $G_1$, extend it by zero to $V(G)$,
and normalize it so that $x_u=\max_{v\in V(G)}x_v=1$. By the eigenvalue equations, we obtain
\begin{equation}\label{eq:eigenvalue-sum}
\begin{aligned}
\rho\sum_{v\in V(G)}x_v
&=\sum_{v\in V(G)}\sum_{w\in N_G(v)}x_w
=\sum_{v\in V(G)}d_G(v)x_v\\
&\le\sum_{v\in V(G)}d_G(v)
=2e(G)\le4n-6.
\end{aligned}
\end{equation}
Using the common-neighbor bound for outerplanar graphs, we obtain
\[
\begin{aligned}
\rho^2=\rho (\rho x_u)
=\rho\sum_{y\in N_G(u)}x_y
&=\sum_{y\in N_G(u)}\sum_{z\in N_G(y)}x_z \\
&=d_G(u)+
  \sum_{z\in V(G)\setminus\{u\}}
  |N_G(u)\cap N_G(z)|x_z \\
&\leq d_G(u)+2\sum_{z\in V(G)\setminus\{u\}}x_z \\
&<d_G(u)+2\sum_{z\in V(G)}x_z \\
&\leq d_G(u)+\frac{4e(G)}{\rho}.
\end{aligned}
\]
By~\eqref{eq:eigenvalue-sum} and $\rho\geq\sqrt{n-1}$, we have
\begin{align*}
d_{G}(u)
&\ge n-1-\frac{8n-12}{\sqrt{n-1}}\\
&=n-1-8\sqrt{n-1}+\frac4{\sqrt{n-1}}\\
&>n-16\sqrt n\ge L.
\end{align*}
The last inequality follows from $n\ge n_0\ge(8L^2+8)^2$ and $L\ge4$.

Choose a vertex $z\notin V(G_1)$.
Delete every edge outside $G_1$ and add $uz$, obtaining an outerplanar graph $G'$.
We claim that $G_1+uz$ is $F$-free.
Otherwise a copy of $F$ must use the new vertex $z$ as a leaf,
with its neighbor mapped to $u$.
Among the $L$ vertices of this copy other than $z$, one is $u$,
so at most $L-1$ neighbors of $u$ are used. Since $d_{G_1}(u)\ge L$,
 an unused neighbor of $u$ in ${G_1}$ can replace $z$.
 We thus obtain a copy of $F$ in ${G_1}$, a contradiction.
  As $F$ is connected and has no isolated vertices,
  $G'$ is also $F$-free. By Lemma~\ref{lem:proper-spanning},
   \[
  \rho(G')\ge\rho(G_1+uz)>\rho(G_1\cup\{z\})=\rho(G_1)=\rho(G),
\]
contrary to the extremality of $G$. Therefore, $G$ is connected.

Since $n-1\ge L$ and $F-w\cong\bigcup_{i=1}^tP_{l_i}$ for the center $w$ of $F$,
we have $F\subseteq K_1\vee P_{n-1}$,
by placing the paths $P_{l_i}$ consecutively along $P_{n-1}$.
Also $F\nsubseteq K_{1,n-1}$ because $F$ is not a star.
Moreover, $64L^4\ge65+L$, so $n\ge(8L^2+8)^2\ge(\sqrt{65+L}+8)^2$.
Since $|V(F)|=L+1$, the order condition in Lemma~\ref{lem:YLstructure} holds.
By Lemma~\ref{lem:YLstructure}, we obtain $G=K_1\vee H$.
\end{proof}

\begin{lemma}\label{lem:edgecomparison}
For an integer $n\ge2$, let $H$ and $H'$ be linear forests
on the same vertex set of order $n-1$.
Let $\rho=\rho(K_1\vee H)$.
If $\rho>3$ and $e(H')-e(H)>4(n-1)/(\rho-3)$, then $\rho(K_1\vee H')>\rho$.
\end{lemma}

\begin{proof}
Let $u$ be a vertex of degree $n-1$ in $K_1\vee H$, and let
$\boldsymbol{x}$ be a positive eigenvector corresponding to
$\rho$, normalized so that $x_u=1$.
For every $v\in V(H)$, by the eigenvalue equation, we obtain
$\rho x_v=x_u +\sum_{w\in N_H(v)}x_w\ge1$.
Dividing by $\rho$, we obtain $x_v\ge1/\rho$.
We also need an upper bound on $x_v$. Choose $v_0\in V(H)$ such that
$x_{v_0}=\max_{v\in V(H)}x_v$.
Since $H$ is a linear forest, $d_H(v_0)\le2$, and each neighbor
$w$ of $v_0$ in $H$ satisfies $x_w\le x_{v_0}$.
Hence $\rho x_{v_0}=x_u +\sum_{w\in N_H(v_0)}x_w
\le1+2x_{v_0}$, which implies
$(\rho-2)x_{v_0}\le1$.
Since $\rho>3$, we have
$x_v\le x_{v_0}\le1/(\rho-2)$ for every $v\in V(H)$.
Therefore,
\[
\frac{1}{\rho}\le x_v\le\frac{1}{\rho-2}
\qquad\text{for every }v\in V(H).
\]

Let $a=|E(H)\setminus E(H')|$ and $b=e(H')-e(H)$.
Then $|E(H')\setminus E(H)|=a+b$. By the Rayleigh principle,
\begin{align*}
\rho(K_1\vee H')-\rho
&\ge\frac{2}{\boldsymbol{x}^{\mathrm T}\boldsymbol{x}}
\left(
\sum_{vw\in E(H')\setminus E(H)}x_vx_w
-\sum_{vw\in E(H)\setminus E(H')}x_vx_w
\right)\\
&\ge\frac{2}{\boldsymbol{x}^{\mathrm T}\boldsymbol{x}}
\left(\frac{a+b}{\rho^2}-\frac{a}{(\rho-2)^2}\right)\\
&=\frac{2}{\rho^2\boldsymbol{x}^{\mathrm T}\boldsymbol{x}}
\left[
b-a\left(\frac{\rho^2}{(\rho-2)^2}-1\right)
\right].
\end{align*}
Since $\rho>3$ and
$(\rho-2)^2>(\rho-1)(\rho-3)$, we obtain
\[
\frac{\rho^2}{(\rho-2)^2}-1
=\frac{4(\rho-1)}{(\rho-2)^2}
<\frac{4}{\rho-3}.
\]
Also, $0\le a\le e(H)\le n-2<n-1$. Thus, by the hypothesis,
\[
b-a\left(\frac{\rho^2}{(\rho-2)^2}-1\right)
\ge b-\frac{4a}{\rho-3}
> b-\frac{4(n-1)}{\rho-3}>0.
\]
Therefore $\rho(K_1\vee H')>\rho$.
\end{proof}

\section{Proof of Theorem~\ref{thm:three}}
We first prove the following lemma.
\begin{lemma}\label{lem:embeddingthree}
For integers $l_1,l_2,l_3$ with $l_1\ge l_2\ge l_3\ge1$,
let $L=l_1+l_2+l_3$.
Let $H=\bigcup_{i=1}^kP_{s_i}$ be a linear forest and let
$n=|V(H)|+1$.
Then $S(l_1,l_2,l_3)\subseteq K_1\vee H$ if and only if at least one of the following conditions holds.
\begin{enumerate}
\item[(i)] $P_L\subseteq P_{s_i}$ for some $i$;
\item[(ii)] for a permutation $(p,q,r)$ of $(l_1,l_2,l_3)$,
there exist distinct $i,j$ such that $P_{p+q}\subseteq P_{s_i}$ and $P_r\subseteq P_{s_j}$;
\item[(iii)] there exist pairwise distinct indices $i_1,i_2,i_3$
such that $P_{l_1}\subseteq P_{s_{i_1}}$, $P_{l_2}\subseteq P_{s_{i_2}}$ and $P_{l_3}\subseteq P_{s_{i_3}}$;
\item[(iv)] for a permutation $(p,q,r)$ of $(l_1,l_2,l_3)$ with $r\ge2$,
there exist distinct $i,j$ such that $P_{p+q+1}\subseteq P_{s_i}$ and $P_{r-1}\subseteq P_{s_j}$.
\end{enumerate}
\end{lemma}
\begin{proof}
We first prove sufficiency. Let $u$ be the vertex of degree $n-1$ corresponding to $K_1$ in $K_1\vee H$.
  Under each of (i)--(iii), we can construct
$S(l_1,l_2,l_3)$ with center $u$.
Under (iv), we can construct it with the center
on $P_{p+q+1}$ and the branch of length $r$
passing through $u$.

We next prove necessity. Let $T$ be a copy of $S(l_1,l_2,l_3)$ with center $w$.
Suppose first that $w=u$.
Then $T-w$ consists of three vertex-disjoint paths $P_{l_1}$, $P_{l_2}$ and $P_{l_3}$ in $H$.
These paths lie in one, two or three components of $H$.
If all three paths lie in the same component $P_{s_i}$,
then $s_i\ge l_1+l_2+l_3=L$, so (i) holds.
If they lie in exactly two components,
then two of them lie in one component and the remaining path lies in the other.
Thus, for a permutation $(p,q,r)$ of $(l_1,l_2,l_3)$,
there exist distinct indices $i,j$ such that $P_{p+q}\subseteq P_{s_i}$ and $P_r\subseteq P_{s_j}$, so (ii) holds.
If they lie in three distinct components,
then there exist pairwise distinct indices $i_1,i_2,i_3$
such that $P_{l_1}\subseteq P_{s_{i_1}}$,
$P_{l_2}\subseteq P_{s_{i_2}}$ and $P_{l_3}\subseteq P_{s_{i_3}}$,
so (iii) holds.

Suppose now that $w\ne u$.
Since $d_H(w)\le2$ and $d_T(w)=3$,
the edge $wu$ belongs to $T$.
Let the path of $T$ beginning with $wu$ have length $r$,
and let the other two paths have lengths $p,q$,
 where $(p,q,r)$ is a permutation of $(l_1,l_2,l_3)$.
 The latter two paths, together with $w$,
 form a path of order $p+q+1$ in some component $P_{s_i}$ of $H$.

If $r=1$, then $p+q+1=L$, so (i) holds.
If $r\ge2$, let $Q$ be the path obtained from this branch by deleting $w$ and $u$.
Then $Q$ has $r-1$ vertices and is vertex-disjoint from the previously identified path on $p+q+1$ vertices. Since $Q$ is contained in $H$, it lies in some component $P_{s_j}$ of $H$.
If $i=j$, then $s_i\ge(p+q+1)+(r-1)=L$, so (i) holds.
If $i\ne j$, then $P_{p+q+1}\subseteq P_{s_i}$ and $P_{r-1}\subseteq P_{s_j}$,
so (iv) holds.
\end{proof}

\begin{proof}[Proof of Theorem~\ref{thm:three}]
Let  $G\in\SPEX(n,S(l_1,l_2,l_3))$.
By Lemma~\ref{lem:reduction},
$G=K_1\vee H$,
where $H=\bigcup_{i=1}^kP_{s_i}$.
Without loss of generality, suppose that
$s_1\ge s_2\ge\cdots\ge s_k\ge1$.
Clearly, $s_i\le L-1$ for every $i$, since otherwise
$S(l_1,l_2,l_3)\subseteq K_1\vee P_{s_i}\subseteq G$,
a contradiction.
Since $n-1\ge(8L^2+8)^2-1=64L^4+128L^2+63>5(L-1)$
and $s_1 \leq L-1$,
we have $k\ge6$.

\claim{1}The paths $P_{s_3},P_{s_4},\ldots,P_{s_{k-1}}$
have the same order, while $P_{s_k}$ may have smaller order. Consequently,
\begin{equation}\label{eq:canonicalthree}
H=P_{s_1}\cup P_{s_2}\cup(k-3)P_{s_3}\cup P_{s_k},
\qquad 1\le s_k\le s_3.
\end{equation}

\begin{proof}
Suppose that two components $P_a$ and $P_b$ of $H$
satisfy $s_3>a\ge b\ge1$.
Apply an $(a,b)$-transformation to $H$,
and denote the resulting forest by $H'$.
Since $a+1\le s_3$, the three largest component orders
of $H'$ are still $s_1,s_2,s_3$.
By Lemma~\ref{lem:embeddingthree},
$K_1\vee H'$ is $S(l_1,l_2,l_3)$-free.
Since $b\le L-1$, we have
\[
n\ge n_0
\ge\max\{1.27\times10^7,6.5025\times2^{L+1}\}
\ge\max\{1.27\times10^7,6.5025\times2^{b+2}\}.
\]
By Lemma~\ref{lem:transfer},
$\rho(K_1\vee H')>\rho(G)$, a contradiction.
Hence at most one component of $H$ has order less than $s_3$.
Thus $H$ has the stated form.
\end{proof}

\claim{2}
\[
s_3=\begin{cases}
l_1-1,& \text{if } l_1\le l_2+l_3+1,\\
l_1-2,& \text{if } l_1\ge l_2+l_3+2.
\end{cases}
\]

\begin{proof}
If $s_3\ge l_1$, Condition (iii) of Lemma~\ref{lem:embeddingthree} holds.
Thus $s_3\le l_1-1$. If $l_1\ge l_2+l_3+2$ and $s_3=l_1-1$,
two distinct paths of order $s_3$ contain $P_{l_2+l_3+1}$ and $P_{l_1-1}$, respectively.
So (iv) holds, a contradiction. Hence $s_3\le l_1-2$.
Let
\[
d=\begin{cases}
l_1-1,& \text{if }l_1\le l_2+l_3+1,\\
l_1-2,& \text{if } l_1\ge l_2+l_3+2.
\end{cases}
\]
Construct a new linear forest $H'=q'P_d\cup P_{r'}$,
where $n-1=q'd+r'$ and $0\le r'<d$. By Lemma~\ref{lem:embeddingthree},
it is easy to verify that $K_1\vee H'$ is $S(l_1,l_2,l_3)$-free.

If $d=1$, then $1\le s_3\le d$ proves the claim.
Otherwise, suppose that $s_3<d$. Let $k$ and $k'$ denote the numbers of components of $H$ and $H'$, respectively.
By~\eqref{eq:canonicalthree} and the construction of $H'$, we obtain
\begin{align*}
k&\ge2+\frac{n-1-s_1-s_2}{s_3}
=\frac{n-1}{s_3}-\frac{s_1+s_2}{s_3}+2
\ge\frac{n-1}{s_3}-2(L-1)+2,
\end{align*}
where $k'\le(n-1)/d+1$. Since $s_3\le d-1$,
\[
\frac1{s_3}-\frac1d\ge\frac1{d-1}-\frac1d
=\frac1{d(d-1)}\ge\frac1{L^2}.
\]
As is well known,
a forest $H$ on $n-1$ vertices with $k$ components has $e(H)=n-1-k$ edges.
Therefore,
\begin{align*}
e(H')-e(H)
&=k-k'\\
&\ge(n-1)\left(\frac1{s_3}-\frac1d\right)-2(L-1)+1\\
&\ge\frac{n-1}{L^2}-2(L-1)+1>\frac{n-1}{L^2}-2(L-1)
>\frac{n-1}{2L^2}.
\end{align*}
For the last inequality, note that
$n-1\ge(8L^2+8)^2-1=64L^4+128L^2+63>2L^2(L^2+1)$, so
$(n-1)/2L^2>L^2+1>2(L-1).$
Since $K_{1,n-1}$ is
$S(l_1,l_2,l_3)$-free, the extremality of $G$
implies $\rho:=\rho(G)\ge\sqrt{n-1}$. Since
\[
n-1-(8L^2+3)^2
\ge(8L^2+8)^2-1-(8L^2+3)^2
=80L^2+54>0,
\]
we have $\rho-3>8L^2$, and therefore
\[
e(H')-e(H)> \frac{n-1}{2L^2} >\frac{4(n-1)}{\rho-3} .
\]
By Lemma~\ref{lem:edgecomparison},
$\rho(K_1\vee H')>\rho(G)$, a contradiction. Thus $s_3=d$, proving Claim~2.
\end{proof}

We now determine $s_1$ and $s_2$.

\case{1} $l_1>l_2$ and $l_1\le l_2+l_3+1$.
Here $s_3=l_1-1\ge l_2\ge l_3$. If $s_1\ge l_1$,
then $P_{l_1}\subseteq P_{s_1}$, $P_{l_2}\subseteq P_{s_2}$ and $P_{l_3}\subseteq P_{s_3}$,
so Condition (iii) of Lemma~\ref{lem:embeddingthree}
implies $S(l_1,l_2,l_3)\subseteq G$,
a contradiction.
Hence $s_1=s_2=s_3=l_1-1$, and by~\eqref{eq:canonicalthree},
we obtain $H=\HH(l_1-1,l_1-1)$. By Lemma~\ref{lem:embeddingthree},
$K_1\vee H$ is $S(l_1,l_2,l_3)$-free.

\case{2} $l_1=l_2+l_3+2$.
Here $s_3=l_1-2=l_2+l_3\ge l_2$.
By Condition (iii) of Lemma~\ref{lem:embeddingthree}, we obtain $s_1\le l_1-1$. If $s_2=l_1-1$,
the first two components both contain $P_{l_1-1}=P_{l_2+l_3+1}$,
contradicting (iv). Hence $s_2=l_1-2$.

Thus $l_1-2\le s_1\le l_1-1$.
Suppose that $s_1=l_1-2$. Apply an $(l_1-2,l_1-2)$-transformation,
and denote the resulting forest by $H_1$.
Then $H_1$ has exactly one component of order $l_1-1$,
while all other components have order at most $l_1-2$.
By Lemma~\ref{lem:embeddingthree}, it is easy to verify that $K_1\vee H_1$ is $S(l_1,l_2,l_3)$-free.
By Lemma~\ref{lem:transfer}, $\rho(K_1\vee H_1)>\rho(G)$, contradicting the extremality of $G$.
Therefore $s_1=l_1-1$, and hence $H=\HH(l_1-1,l_1-2)$.

\case{3} $l_1\ge l_2+l_3+3$.
Here $s_3=l_1-2\ge l_2$.
Thus Condition (iii) of Lemma~\ref{lem:embeddingthree} implies $s_1\le l_1-1$.
If $s_1=l_1-1$, then $P_{l_1-1}\subseteq P_{s_1}$ and
$ P_{l_2+l_3+1}\subseteq P_{s_3}$, contradicting (iv).
Thus $s_1=s_2=s_3=l_1-2$, and we obtain $H=\HH(l_1-2,l_1-2)$.
By Lemma~\ref{lem:embeddingthree}, $K_1 \vee H$ is $S(l_1,l_2,l_3)$-free.

\case{4} $l_1=l_2>l_3$.
Now $s_3=l_1-1$. If $s_2\ge l_1$,
the first two components contain disjoint $P_{l_1}$ and the third contains $P_{l_3}$,
contradicting (iii).
Thus $s_2=l_1-1$. If $s_1\ge l_1+l_3+1$,
then the first and second components contain $P_{l_1+l_3+1}$ and $P_{l_1-1}$,
contradicting (iv).
Hence $s_1\le l_1+l_3$.

Suppose that $s_1<l_1+l_3$.
Apply an $(s_1,l_1-1)$-transformation, using one of the repeated paths in~\eqref{eq:canonicalthree},
and denote the resulting forest by $H_2$.
Then $H_2$ has at most one component of order at least $l_1$,
and every component has order at most $l_1+l_3$. By Lemma~\ref{lem:embeddingthree},
it is easy to verify that $K_1\vee H_2$ is $S(l_1,l_2,l_3)$-free.
By Lemma~\ref{lem:transfer}, $\rho(K_1\vee H_2)>\rho(G)$, contradicting the extremality of $G$.
Therefore $s_1=l_1+l_3$, and hence $H=\HH(l_1+l_3,l_1-1)$. Again by Lemma~\ref{lem:embeddingthree},
$K_1\vee H$ is $S(l_1,l_2,l_3)$-free.

In each case, $K_1\vee H^*$ is an
$S(l_1,l_2,l_3)$-free outerplanar graph,
and every extremal graph is isomorphic to $K_1\vee H^*$.
\end{proof}

\section{Proof of Theorem~\ref{thm:many}}
\begin{proof}[Proof of Theorem~\ref{thm:many}]
Let $G\in\SPEX(n,S(l_1,l_2,\ldots,l_t))$, where $t\geq 4$.
By Lemma~\ref{lem:reduction},
$G=K_1\vee H$,
where $H=\bigcup_{i=1}^kP_{s_i}$.
Without loss of generality, suppose that
$s_1\ge s_2\ge\cdots\ge s_k\ge1$.
Clearly, $s_i\le L-1$ for every $i$, since otherwise
$S(l_1,l_2,\ldots,l_t)\subseteq K_1\vee P_{s_i}\subseteq G$,
a contradiction.
For any linear forest $J$, every vertex of $J$
has degree at most $3$ in $K_1\vee J$.
Since $t\ge4$, the center of any starlike tree
$S(l_1,l_2,\ldots,l_t)$ contained in $K_1\vee J$
must be the vertex corresponding to $K_1$.
Deleting this vertex leaves
$\bigcup_{i=1}^tP_{l_i}$ in $J$.
Conversely, if these paths are contained in $J$,
joining one endpoint of each path to the vertex
corresponding to $K_1$ yields $S(l_1,l_2,\ldots,l_t)$.
Therefore, for every linear forest $H$, we have
\begin{equation}\label{eq:manyembedding}
S(l_1,l_2,\ldots,l_t)\subseteq K_1\vee H
\quad\Longleftrightarrow\quad
\bigcup_{i=1}^tP_{l_i}\subseteq H.
\end{equation}
Let
$H_0=H_{\mathcal{OP}}(hl_1-1,l_1-1).$
By its definition, $H_0$ has one path component of order $hl_1-1$,
while every other component has order at most $l_1-1$. Hence, every
copy of $P_{l_1}$ in $H_0$ must lie in the path component of order
$hl_1-1$. This component contains at most
$
\left\lfloor\frac{hl_1-1}{l_1}\right\rfloor=h-1
$
pairwise vertex-disjoint copies of $P_{l_1}$.
By the definition of $h$,
the linear forest $\bigcup_{i=1}^{t}P_{l_i}$ contains $h$ vertex-disjoint copies of $P_{l_1}$. Therefore,
$\bigcup_{i=1}^{t}P_{l_i}\not\subseteq H_0,$
and hence $K_1\vee H_0$ is $S(l_1,l_2,\ldots,l_t)$-free by~\eqref{eq:manyembedding}.

\claim{1} $H$ has at least $t$ path components of order at least $l_1-1$.

\begin{proof}
Let $k$ denote the number of components of $H$.
Since every component of $H$ has order at most $L-1$,
    we have $n-1=\sum_{i=1}^k s_i\le k(L-1)$.
    Since $l_i\ge1$ for $1\le i\le t$, we have
$L=\sum_{i=1}^t l_i\ge t$.
Together with $n\ge n_0$, this yields
$n-1\ge(8L^2+8)^2-1
>L(L-1)\ge t(L-1).$
Hence $k\ge(n-1)/(L-1)>t$.
    If $l_1=2$, then every component of $H$ has order at least $l_1-1=1$.
Thus the claim holds when $l_1=2$.

For $l_1\ge3$,
suppose to the contrary that
at most $t-1$ components of $H$ have order
at least $l_1-1$.
    The $t-1$ largest components each have order at most $L-1$,
    while each of the remaining $k-t+1$ components
    has order at most $l_1-2$.
    Hence
\[
\begin{aligned}
n-1
&\le(t-1)(L-1)+(k-t+1)(l_1-2)\\
&\le(t-1)(L-1)+k(l_1-2),
\end{aligned}
\]
which implies
\[
k\ge\frac{n-1-(t-1)(L-1)}{l_1-2}.
\]

Let $k^*$ be the number of components of $H_0$.
By its definition, we have $k^*\le\frac{n-1}{l_1-1}+1$.
Since $t\le L$, $l_1-2\ge1$ and $(l_1-1)(l_1-2)\le L^2$, we obtain
\begingroup
\allowdisplaybreaks[4]
\begin{align*}
e(H_0)-e(H)
&=k-k^*\\
&\ge\frac{n-1-(t-1)(L-1)}{l_1-2}
-\frac{n-1}{l_1-1}-1\\
&=\frac{n-1}{(l_1-1)(l_1-2)}
-\frac{(t-1)(L-1)}{l_1-2}-1\\
&\ge\frac{n-1}{L^2}-(t-1)(L-1)-1\\
&\ge\frac{n-1}{L^2}-(L-1)^2-1\\
&\ge\frac{n-1}{L^2}-L^2-1\\
&>\frac{n-1}{L^2}-\frac{n-1}{2L^2}
=\frac{n-1}{2L^2}.
\end{align*}
\endgroup
The strict inequality holds because $n-1\ge(8L^2+8)^2-1=64L^4+128L^2+63>2L^2(L^2+1)$. Furthermore,
$
n-1-(8L^2+3)^2
\ge(8L^2+8)^2-1-(8L^2+3)^2
=80L^2+54>0.
$
Since $K_{1,n-1}$ is
$S(l_1,l_2,\ldots,l_t)$-free, the extremality
of $G$ implies $\rho(G)\ge\sqrt{n-1}$.
Hence $\rho(G)-3>8L^2$. Thus
\[
e(H_0)-e(H)>\frac{n-1}{2L^2}=\frac{4(n-1)}{8L^2}>
\frac{4(n-1)}{\rho(G)-3}.
\]
By Lemma~\ref{lem:edgecomparison},
we obtain $\rho(K_1\vee H_0)>\rho(G)$, a contradiction. This proves Claim~1.
\end{proof}

Recall that $l_i=l_1$ for $1\le i\le h$
    and $l_i\le l_1-1$ for  $h<i\le t$.
    We now give a criterion for $K_1\vee J$
    to be $S(l_1,l_2,\ldots,l_t)$-free.
\claim{2} Let $J=\bigcup_{i=1}^{k'}P_{r_i}$ be a linear forest in which
at least $t$ components are of order at least $l_1-1$. Then
\begin{equation}\label{eq:capacity}
K_1\vee J\text{ is }S(l_1,l_2,\ldots,l_t)\text{-free}
\quad\Longleftrightarrow\quad
\sum_{i=1}^{k'}\floor{\frac{r_i}{l_1}}\le h-1.
\end{equation}
\begin{proof}
The maximum number of pairwise vertex-disjoint subpaths
    of order $l_1$ in $J$ is
    $\sum_{i=1}^{k'}\lfloor r_i/l_1\rfloor$.
If this sum is at most $h-1$, then
$\bigcup_{i=1}^tP_{l_i}\nsubseteq J$,
since $h$ of these paths have order $l_1$.
By~\eqref{eq:manyembedding},
$K_1\vee J$ is $S(l_1,l_2,\ldots,l_t)$-free.

Conversely, if the sum is at least $h$, choose $h$
pairwise vertex-disjoint subpaths of order $l_1$ in $J$.
These subpaths lie in at most $h$ components of $J$.
Hence at least $t-h$ components of order at least $l_1-1$
remain unused.
Since $l_i\le l_1-1$ for $h<i\le t$, each remaining $P_{l_i}$
can be chosen in a distinct unused component.
Therefore $\bigcup_{i=1}^tP_{l_i}\subseteq J$,
and~\eqref{eq:manyembedding} implies
$S(l_1,l_2,\ldots,l_t)\subseteq K_1\vee J$.
\end{proof}

\claim{3} At most one path component of $H$ has order at least $l_1$.

\begin{proof}
Suppose that $H$ has two components $P_a$ and $P_b$
with $a\ge b\ge l_1$.
Repeatedly apply the transformation, transferring one vertex
from the second path to the first until the second path
has order $l_1-1$.
Denote the resulting forest by $H'$.
Then $P_a\cup P_b$ is replaced by
$P_{a+b-l_1+1}\cup P_{l_1-1}$.
Indeed, the second path loses exactly $b-(l_1-1)$ vertices,
all of which are added to the first.
Since $b\le L-1$, we have
\(n\ge n_0
\ge\max\{1.27\times10^7,6.5025\times2^{L+1}\}
\ge\max\{1.27\times10^7,6.5025\times2^{b+2}\}.\)
 By Lemma~\ref{lem:transfer}, we obtain
 $\rho(K_1\vee H')>\rho(G)$.
Both new paths have order at least $l_1-1$,
so $H'$ still has at least $t$ components
of order at least $l_1-1$.
  Thus $H'$ satisfies the hypothesis of Claim~2.

Let $\alpha,\beta,\gamma,\delta$ be integers such that $a=\alpha l_1+\beta$ and $b=\gamma l_1+\delta$,
where $0\le\beta,\delta\le l_1-1$. Since
$\beta+\delta-l_1+1\le2(l_1-1)-l_1+1=l_1-1,$
we have
\begin{align*}
\floor{\frac{a+b-l_1+1}{l_1}}
&=\alpha+\gamma+\floor{\frac{\beta+\delta-l_1+1}{l_1}}\\
&\le\alpha+\gamma+\floor{\frac{l_1-1}{l_1}}=\alpha+\gamma
=\floor{\frac{a}{l_1}}+\floor{\frac{b}{l_1}}.
\end{align*}
The other new path contributes $\floor{(l_1-1)/l_1}=0$.
Hence the sum in~\eqref{eq:capacity} does not increase.
By Claim~2, $K_1\vee H'$ remains
$S(l_1,l_2,\ldots,l_t)$-free,
contradicting the extremality of $G$.
This proves Claim~3.
\end{proof}

\claim{4} $s_1=hl_1-1$, and every other component has order $l_1-1$,
except possibly one of smaller order.

\begin{proof}
By Claims~2 and~3, we have $s_1\le hl_1-1$,
and every other component has order at most $l_1-1$.
If $s_1<hl_1-1$, another path $P_b$ exists, since $s_1\le L-1<n-1$.
Apply an $(s_1,b)$-transformation and denote the resulting forest by $H'$. The first path still has order at most $hl_1-1$,
and all other paths have order at most $l_1-1$.
Thus $H'$  contains at most $h-1$
pairwise vertex-disjoint subpaths of order $l_1$.
By~\eqref{eq:manyembedding}, $K_1\vee H'$ is $S(l_1,l_2,\ldots,l_t)$-free.
By Lemma~\ref{lem:transfer},
the spectral radius strictly increases, contradicting extremality.
Hence $s_1=hl_1-1$.

If two components other than $P_{s_1}$ have orders $a,b<l_1-1$,
with $a\ge b\ge1$, apply an $(a,b)$-transformation.
By Claim~2, the resulting graph remains $S(l_1,l_2,\ldots,l_t)$-free,
while its spectral radius strictly increases by Lemma~\ref{lem:transfer},
a contradiction.
Thus at most one component has order less than $l_1-1$, proving Claim~4.
\end{proof}

It follows that $H=P_{hl_1-1}\cup qP_{l_1-1}\cup P_r$, where $n-hl_1=q(l_1-1)+r$ and $0\le r<l_1-1$. Therefore $G\cong K_1\vee\HH(hl_1-1,l_1-1)$.
\end{proof}

\section{Proofs of Theorems~\ref{thm:bouquet} and~\ref{thm:forest}}
We now prove the result for $C_{l_1,l_2,\ldots,l_t}$ stated in the introduction.
\begin{proof}[Proof of Theorem~\ref{thm:bouquet}]
Let $B=C_{l_1,l_2,\ldots,l_t}$,
$F_0=\bigcup_{i=1}^tP_{l_i-1}$, and
$T=S(l_1-1,l_2-1,\ldots,l_t-1,1,1)$.
The sum of the branch lengths of $T$ is
$L=\sum_{i=1}^t(l_i-1)+2$.
Thus $n\ge n_0$ ensures that the order conditions
in Lemma~\ref{lem:YLstructure} and
Theorem~\ref{thm:many} hold, and
$n-1\ge |V(F_0)|+2$.

Let $G\in\SPEX(n,B)$.
We first show that $G$ is connected.
Otherwise, let $G_1$ be a component with
$\rho(G_1)=\rho(G)$.
Choose $u\in V(G_1)$ and $z\in V(G)\setminus V(G_1)$.
Delete all edges outside $G_1$ and add $uz$,
obtaining an outerplanar graph $G'$.
Since every vertex of $B$ has degree at least $2$,
any subgraph isomorphic to $B$ in $G'$ would lie in $G_1$.
Hence $G'$ is $B$-free.
By Lemma~\ref{lem:proper-spanning},
\[
\rho(G')=\rho(G_1+uz)>\rho(G_1)=\rho(G),
\]
a contradiction.

Since $B\subseteq K_1\vee P_{n-1}$ and
$B\nsubseteq K_{1,n-1}$,
Lemma~\ref{lem:YLstructure} with $t=1$ implies
$G=K_1\vee H$, where $H$ is a linear forest
of order $n-1$.

For any linear forest $H$ of order $n-1$,
let $u$ be the vertex corresponding to $K_1$.
If $B\subseteq K_1\vee H$, its vertex of degree $2t$
must be $u$, since every vertex of $H$ has degree
at most $3$ in $K_1\vee H$.
Deleting $u$ from $B$ leaves $F_0$, so $F_0\subseteq H$.
Conversely, if $F_0\subseteq H$, joining the two endpoints
of each path to $u$ yields $B$.
Therefore,
$
B\subseteq K_1\vee H$ if and only if $F_0\subseteq H.$

Similarly, since $T$ has $t+2\ge4$ branches, its center
    must be $u$ whenever $T\subseteq K_1\vee H$.
    Hence $T\subseteq K_1\vee H$ if and only if
    $F_0\cup2P_1\subseteq H$.
    If $F_0\subseteq H$, then at least two vertices of $H$
    lie outside this subgraph, since $n-1\ge |V(F_0)|+2$.
    These two vertices can be chosen as the additional $2P_1$,
    so $F_0\cup2P_1\subseteq H$.
    Thus,
    $F_0\subseteq H$ if and only if $F_0\cup2P_1\subseteq H$.
    Consequently,
\[
K_1\vee H \text{ is } B\text{-free}
\quad\Longleftrightarrow\quad
K_1\vee H \text{ is } T\text{-free}.
\]
The maximum branch length of $T$ is $l_1-1$
with multiplicity $h$.
By Theorem~\ref{thm:many}, the unique $T$-free
extremal graph is
$G^*=K_1\vee\HH(h(l_1-1)-1,l_1-2)$.
Since $G$ is $B$-free, it is also $T$-free.
By Theorem~\ref{thm:many}, we have
$\rho(G)\le\rho(G^*)$.
Conversely, $G^*$ is $T$-free and hence $B$-free.
The extremality of $G$ therefore implies
$\rho(G^*)\le\rho(G)$.
Thus $\rho(G)=\rho(G^*)$, and the equality condition
in Theorem~\ref{thm:many} yields $G\cong G^*$.
\end{proof}

For sufficiently large $n$, Theorem~\ref{thm:bouquet} yields the same
extremal graphs as those in Theorems~\ref{thm:knowncycles}
and~\ref{thm:twocycles}, and completely resolves
Problem~\ref{prob:cycles} for outerplanar graphs.
In particular, for $C_{l_1,l_2,l_3}$ with
$l_1\ge l_2\ge l_3\ge3$, the unique extremal graph is
$K_1\vee H^*$ whenever $n$ satisfies the bound in
Theorem~\ref{thm:bouquet}.
The corresponding forests $H^*$ are listed
in Table~\ref{tab:bouquet}.
\begin{table}[ht]
\centering
\caption{Extremal forests for $C_{l_1,l_2,l_3}$.}\label{tab:bouquet}\smallskip
\begin{tabular}{@{}ll@{}}
\toprule
Cycle lengths & $H^*$\\
\midrule
$l_1>l_2$ & $\HH(l_1-2,l_1-2)$\\
$l_1=l_2>l_3$ & $\HH(2l_1-3,l_1-2)$\\
$l_1=l_2=l_3$ & $\HH(3l_1-4,l_1-2)$\\
\bottomrule
\end{tabular}
\end{table}

The starlike-tree results can also be applied to certain
linear forests. We now prove the corresponding result
stated in the introduction.

\begin{proof}[Proof of Theorem~\ref{thm:forest}]
Let $F=F(l_1,l_2,\ldots,l_t)$,
and let $G\in\SPEX(n,F)$.
In each case, we choose a starlike tree $T$ containing $F$.
The sum of its branch lengths is at most $L+1$.
Thus the bound on $n$ in Theorem~\ref{thm:forest} is at least
$n_0$ in~\eqref{eq:n0} for each such tree.
Since every $F$-free graph is $T$-free, we have
$\operatorname{spex}_{\mathcal{OP}}(n,F)
\le \operatorname{spex}_{\mathcal{OP}}(n,T)$.
Let $G^*=K_1\vee J$ be the unique
$n$-vertex $T$-free outerplanar graph
with maximum spectral radius.
Let $u$ be the vertex corresponding to $K_1$.
If $G^*$ is $F$-free, then
$\operatorname{spex}_{\mathcal{OP}}(n,F)
\ge \rho(G^*)
=\operatorname{spex}_{\mathcal{OP}}(n,T)$,
so equality holds.
Every $F$-free extremal graph is then also a $T$-free
extremal graph and hence is isomorphic to $G^*$.
Thus, in each case, it suffices to choose $T$ containing $F$
and verify that the unique $T$-free extremal graph
is $F$-free.

\textbf{Case 1.} $l_1=l_2$.
Let $T$ be the starlike tree with branch lengths
$l_i$ for $2\le i\le t$, together with $l_1-1,1,1$.
Clearly, $F\subseteq T$.
The maximum branch length of $T$ is $l_1$
with multiplicity $h-1$.
Since $T$ has $t+2\ge4$ branches,
Theorem~\ref{thm:many} implies
$J=\HH((h-1)l_1-1,l_1-1)$.
Thus $J$ contains at most $h-2$ pairwise
vertex-disjoint subpaths of order $l_1$.

Suppose that $F\subseteq G^*$.
Among the $h$ paths of order $l_1$ in $F$,
at most one contains $u$.
Hence at least $h-1$ of these paths lie in
$J=G^*-u$ and are pairwise vertex-disjoint,
a contradiction.
Therefore, $G^*$ is $F$-free.

\textbf{Case 2.} $l_2<l_1\le2l_2-1$.
Let $T$ be the starlike tree with branch lengths
$l_i$ for $2\le i\le t$, together with
$l_2-1$ and $l_1-l_2$.
The last two branches, together with the center,
form $P_{l_1}$, since
$(l_2-1)+(l_1-l_2)+1=l_1$.
Deleting the edges from the center to the other
branches leaves $P_{l_2},P_{l_3},\ldots,P_{l_t}$.
Thus $F\subseteq T$.

Recall that $h$ is the number of indices
$i\in\{2,3,\ldots,t\}$ for which $l_i=l_2$.
If $t\ge3$, then $T$ has $t+1\ge4$ branches,
and its maximum branch length is $l_2$
with multiplicity $h$.
By Theorem~\ref{thm:many},
$J=\HH(hl_2-1,l_2-1)$.
If $t=2$, then $h=1$ and
$T=S(l_2,l_2-1,l_1-l_2)$.
Both $l_2-1$ and $l_1-l_2$ are less than $l_2$,
and $l_2<(l_2-1)+(l_1-l_2)+1=l_1$.
By the first case of Theorem~\ref{thm:three},
we obtain the same forest $J=\HH(hl_2-1,l_2-1)$.

Let $Q$ be the component of $J$ of order $hl_2-1$.
Every other component has order at most $l_2-1$.
Thus every subpath of $J$ of order at least $l_2$
lies in $Q$.
Suppose that $F\subseteq G^*$, and consider its
pairwise vertex-disjoint paths
$P_{l_1},P_{l_2},\ldots,P_{l_t}$.
If $u\in V(P_{l_1})$, then all $h$ paths of order
$l_2$ lie in $Q$.
Hence $|V(Q)|\ge hl_2$, contradicting
$|V(Q)|=hl_2-1$.
If $u\notin V(P_{l_1})$, then $P_{l_1}$ and at least
$h-1$ paths of order $l_2$ lie in $Q$.
Therefore,
$|V(Q)|\ge l_1+(h-1)l_2>hl_2-1$,
a contradiction.
Thus $G^*$ is $F$-free.

\textbf{Case 3.} $l_1=2l_2$.
Let $T$ be the starlike tree with branch lengths
$l_i$ for $2\le i\le t$, together with
$l_2$ and $l_2-1$.
Similarly, the last two branches and the center
form $P_{l_1}$, so $F\subseteq T$.
If $t\ge3$, then $T$ has $t+1\ge4$ branches,
and its maximum branch length is $l_2$
with multiplicity $h+1$.
By Theorem~\ref{thm:many},
$J=\HH((h+1)l_2-1,l_2-1)$.
If $t=2$, then $h=1$ and
$T=S(l_2,l_2,l_2-1)$.
By the fourth case of Theorem~\ref{thm:three},
we obtain the same forest
$J=\HH((h+1)l_2-1,l_2-1)$.

Let $Q$ be the component of $J$ of order $(h+1)l_2-1$.
Every other component has order at most $l_2-1$,
so every subpath of $J$ of order at least $l_2$
lies in $Q$.
Suppose that $F\subseteq G^*$, and consider its
pairwise vertex-disjoint paths
$P_{l_1},P_{l_2},\ldots,P_{l_t}$.
If $u\in V(P_{l_1})$, then deleting $u$ from
$P_{l_1}$ leaves at most two paths containing
$2l_2-1$ vertices in total.
At least one has order at least $l_2$ and hence
lies in $Q$.
All $h$ paths of order $l_2$ also lie in $Q$
and are vertex-disjoint from this path.
Thus $|V(Q)|\ge l_2+hl_2=(h+1)l_2$,
a contradiction.
If $u\notin V(P_{l_1})$, then $P_{l_1}$ and at least
$h-1$ paths of order $l_2$ lie in $Q$.
Hence $|V(Q)|\ge l_1+(h-1)l_2=(h+1)l_2$,
 a contradiction.
Therefore, $G^*$ is $F$-free.

\textbf{Case 4.} $l_1>2l_2$ and $l_1$ is even.
Let $T$ be the starlike tree with branch lengths
$l_i$ for $2\le i\le t$, together with
$l_1/2$ and $l_1/2-1$.
Similarly, the last two branches and the center
form $P_{l_1}$, so $F\subseteq T$.
Since $l_1/2>l_2$, the maximum branch length of $T$
is $l_1/2$ with multiplicity one.
If $t\ge3$, then $T$ has $t+1\ge4$ branches.
By Theorem~\ref{thm:many},
$J=\HH(l_1/2-1,l_1/2-1)$.
If $t=2$, then
$T=S(l_1/2,l_1/2-1,l_2)$.
Its maximum branch length is unique and satisfies
$l_1/2<(l_1/2-1)+l_2+1$.
By the first case of Theorem~\ref{thm:three},
we obtain the same forest
$J=\HH(l_1/2-1,l_1/2-1)$.

Every component of $J$ has order at most $l_1/2-1$.
We show that $P_{l_1}\nsubseteq G^*$,
which implies that $G^*$ is $F$-free.
Suppose that $P_{l_1}\subseteq G^*$.
If this path avoids $u$, it lies in a component of $J$,
contradicting the bound on its order.
If it contains $u$, deleting $u$ leaves at most two
paths in $J$, each of order at most $l_1/2-1$.
This implies $l_1-1\le2(l_1/2-1)=l_1-2$,
a contradiction.
Therefore, $G^*$ is $F$-free.

In each case, $F\subseteq T$ and $G^*$ is $F$-free.
By the argument above, $G^*$ is the unique
$F$-free extremal graph.
This completes the proof.
\end{proof}
For sufficiently large $n$, Theorem~\ref{thm:forest} yields the same extremal graphs as those determined by Yin, Li and Meng~\cite{YLM} for $tP_l$ with $t\ge2$ and $l\ge3$, and extends their result to linear forests whose path orders are not all equal.
For $t\ge2$ and $l_t\ge2$,
the case $l_1>2l_2$ with $l_1$ odd is not covered by Theorem~\ref{thm:forest}.
The constructions used above do not settle this case through Theorems~\ref{thm:three} and~\ref{thm:many}.
For example, $F(7,2,2,2)\subseteq S(5,2,2,2,1)$,
but by Theorem~\ref{thm:many}, the unique extremal graph
for this tree is $K_1\vee H_{\mathcal{OP}}(4,4)$,
which contains $F(7,2,2,2)$ for sufficiently large $n$.
Thus this choice of tree does not yield the extremal graph
for $F(7,2,2,2)$-free outerplanar graphs.

A theta graph $\theta_{r,p,q}$ consists of three internally vertex-disjoint paths
with the same endpoints. If all three paths have length at least $2$,
it contains a subdivision of $K_{2,3}$ and is therefore not outerplanar.
Thus only $\theta_{1,p,q}$ with $q\ge p\ge2$ can occur in an
outerplanar graph. The following corollary consequently treats all
theta graphs that can occur in outerplanar graphs.

Wang, Shi and Sun~\cite[Theorem~1.9]{WSS} characterized
$K_1\vee\HH(p+q-2,p+q-2)$ as the unique extremal graph
for $\theta_{1,p,q}$-free outerplanar graphs whenever
$n\ge\max\{400,5\cdot2^{p+q-1}+1\}$.
The following corollary gives an alternative derivation of their extremal graph using Theorem~\ref{thm:many}.
\begin{corollary}\label{cor:theta}
For integers $p,q$ with $q\ge p\ge2$, let $L=p+q+2$.
For an integer $n\ge n_0$ \textup{(see~\eqref{eq:n0})}, let $G$ be an $n$-vertex
$\theta_{1,p,q}$-free outerplanar graph. Then
$\rho(G)\le\rho(K_1\vee\HH(p+q-2,p+q-2))$,
with equality if and only if
$G\cong K_1\vee\HH(p+q-2,p+q-2)$.
\end{corollary}
\begin{proof}
Let $T=S(p+q-1,1,1,1)$.
Let $G$ be an $n$-vertex $\theta_{1,p,q}$-free
outerplanar graph with maximum spectral radius.
A similar argument shows that $G$ is connected.
Since $\theta_{1,p,q}$ is a subgraph of
$K_1\vee P_{n-1}$ but is not a subgraph of
$K_{1,n-1}$, by Lemma~\ref{lem:YLstructure},
we have $G=K_1\vee H$, where $H$ is a linear forest.
Let $u$ be the vertex of degree $n-1$.

We first show that $G$ contains $\theta_{1,p,q}$
if and only if $H$ contains $P_{p+q-1}$.
Suppose that $\theta_{1,p,q}\subseteq G$.
This theta subgraph must contain $u$, since $H$
contains no cycle.
Deleting any vertex from $\theta_{1,p,q}$ leaves
a connected graph.
Its remaining $p+q-1$ vertices therefore lie
in one path component of $H$, so
$P_{p+q-1}\subseteq H$.
Conversely, suppose that $H$ contains the path
$v_1v_2\cdots v_{p+q-1}$.
Then the edge $uv_p$ and the two paths
$uv_1v_2\cdots v_p$ and
$uv_{p+q-1}v_{p+q-2}\cdots v_p$
form $\theta_{1,p,q}$.

We also have $T\subseteq G$ if and only if
$P_{p+q-1}\subseteq H$.
Indeed, every vertex of $H$ has degree at most $3$
in $G$, whereas the center of $T$ has degree $4$.
Thus the center must be $u$, and deleting it
from the longest branch leaves $P_{p+q-1}$ in $H$.
Conversely, given a path $P_{p+q-1}$ in $H$,
there are at least three vertices outside this path
in $H$, since $n\ge p+q+3$.
Together with $u$, these vertices and the path
contain $T$.

Let $G^*=K_1\vee\HH(p+q-2,p+q-2)$.
Since $G$ is $T$-free, Theorem~\ref{thm:many}
implies that $\rho(G)\leq\rho(G^{*})$,
with equality only if $G\cong G^{*}$.
Conversely, every component of
$\HH(p+q-2,p+q-2)$ has order at most $p+q-2$,
so it contains no $P_{p+q-1}$.
By the argument above, $G^{*}$ is
$\theta_{1,p,q}$-free.
By the choice of $G$, we have
$\rho(G)\geq\rho(G^{*})$.
Therefore, $\rho(G)=\rho(G^{*})$, and hence
$G\cong G^{*}$.
\end{proof}

Yin and Li~\cite[Theorem~3]{YL} and
Wang, Shi and Sun~\cite[Theorem~1.6]{WSS} determined the same
unique extremal graphs for outerplanar graphs forbidding
$(k+1)K_2$. Wang, Shi and Sun proved their result for
$n\ge\max\{400,5\cdot2^{2k}+1\}$.
The following corollary gives an alternative derivation of these
extremal graphs using Theorems~\ref{thm:three} and~\ref{thm:many},
with a larger sufficient bound on $n$ than that of Wang, Shi and Sun.
\begin{corollary}\label{cor:matching}
For an integer $k\ge1$, let $L=2k+2$.
For an integer $n\ge n_0$ \textup{(see~\eqref{eq:n0})}, let $G$ be an $n$-vertex
$(k+1)K_2$-free outerplanar graph. Then
$\rho(G)\le\rho(K_1\vee\HH(2k-1,1))$,
with equality if and only if $G\cong K_1\vee\HH(2k-1,1)$.
\end{corollary}
\begin{proof}
Let $T=S(\underbrace{2,2,\ldots,2}_{k},1,1)$.
The tree $T$ contains $(k+1)K_2$.
To see this, select the edge not incident with the
center from each branch of length $2$,
together with one branch of length $1$.
These $k+1$ edges have pairwise disjoint endpoints.
Thus every $(k+1)K_2$-free graph is $T$-free.
Let $G$ be an $n$-vertex $(k+1)K_2$-free
outerplanar graph with maximum spectral radius.

If $k=1$, then $T=S(2,1,1)$.
By Theorem~\ref{thm:three}, the unique $T$-free
outerplanar extremal graph is $K_{1,n-1}$.
Since $K_{1,n-1}$ is also $2K_2$-free,
it is the unique $2K_2$-free outerplanar
extremal graph, as required.

It remains to consider $k\ge2$. Then $T$ has $k+2$ branches,
exactly $k$ of which have maximum length $2$.
By Theorem~\ref{thm:many}, its unique outerplanar
extremal graph is
$G^*=K_1\vee\HH(2k-1,1)
=K_1\vee\bigl(P_{2k-1}\cup(n-2k)K_1\bigr)$.
Let $u$ be the vertex of degree $n-1$ in $G^*$.
Any matching in $G^*$ contains at most one edge
incident with $u$.
All remaining edges of the matching lie in
$P_{2k-1}$, which contains at most $k-1$
pairwise vertex-disjoint edges.
Hence $G^*$ is $(k+1)K_2$-free.

Since every $(k+1)K_2$-free graph is $T$-free,
Theorem~\ref{thm:many} implies that
$\rho(G)\leq\rho(G^{*})$.
Conversely, $G^{*}$ is $(k+1)K_2$-free.
By the choice of $G$, we have $\rho(G)\geq\rho(G^{*})$.
Therefore, $\rho(G)=\rho(G^{*})$, and hence
$G\cong G^{*}$.
\end{proof}

\section*{Data availability}

No data were used for the research described in this article.

\section*{Declaration of competing interest}

The authors declare that they have no conflict of interest.

\end{document}